\documentclass[12pt]{amsart}
\newtheorem{thm}{Theorem}[section]
\newtheorem{theorem}[thm]{Theorem}
\newtheorem{corollary}[thm]{Corollary}
\newtheorem{lemma}[thm]{Lemma}
\newtheorem{proposition}[thm]{Proposition}

\theoremstyle{definition}

\theoremstyle{remark}

\numberwithin{equation}{section}

\newcommand{\ov}{\overline}
\newcommand{\wt}{\widetilde}
\newcommand{\bn}{\binom}

\newcommand{\bC}{{\mathbb C}}
\newcommand{\bN}{{\mathbb N}}
\newcommand{\bR}{{\mathbb R}}

\newcommand{\cH}{{\mathcal H}}

\newcommand{\fD}{{\mathfrak D}}

\newcommand{\al}{\alpha}

\newcommand{\ga}{\gamma}
\newcommand{\de}{\delta}
\newcommand{\la}{\lambda}
\newcommand{\si}{\sigma}
\renewcommand{\th}{\theta}

\begin{document}

\title[Reconstruction by two spectra]%
{Inverse spectral problems for Sturm-Liouville operators with
singular potentials, II.\\ Reconstruction by two spectra${}^{\dag}$}%
\author[R.~O.~Hryniv, Ya.~V.~Mykytyuk]{Rostyslav O.~Hryniv and Yaroslav V.~Mykytyuk}%
\address{Institute for Applied Problems of Mechanics and Mathematics,
3b~Naukova st., 79601 Lviv, Ukraine and Lviv National University,
1 Universytetska st.,
79602 Lviv, Ukraine}%
\email{rhryniv@iapmm.lviv.ua}%

\address{Lviv National University, 1 Universytetska st., 79602 Lviv, Ukraine}%
\email{yamykytyuk@yahoo.com}%

\thanks{${}^{\dag}$The work was partially supported by Ukrainian Foundation
for Basic Research DFFD under  grant No.~01.07/00172}%

\date{\today}%

\begin{abstract}
We solve the inverse spectral problem of recovering the singular potentials
$q\in W^{-1}_{2}(0,1)$ of Sturm-Liouville operators by two spectra. The
reconstruction algorithm is presented and necessary and sufficient
conditions on two sequences to be spectral data for Sturm-Liouville
operators under consideration are given.
\end{abstract}
\maketitle


\section{Introduction}

Suppose that $q$ is a real-valued distribution from the space
$W^{-1}_2(0,1)$. A Sturm-Liouville operator $T$ in a Hilbert space $\cH :=
L_2(0,1)$ corresponding formally to the differential expression
\begin{equation}\label{eq:de}
    -\frac{d^2}{dx^2} + q
\end{equation}
can be defined as follows~\cite{SS}. We take a distributional primitive
$\si\in \cH$ of~$q$, put
\[
    l_\si(u) := - (u'-\si u)' - \si u'
\]
for absolutely continuous functions $u$ whose {\em quasi-derivative}
$u^{[1]}:= u' -\si u$ is again absolutely continuous, and then set
\[
    T u := l_\si (u)
\]
for suitable $u$ (see Section~\ref{sec:asymp} for rigorous definitions).
When considered on this natural domain subject to the boundary conditions
\begin{equation}\label{eq:bc}
    u^{[1]}(0) - H u(0) = 0, \quad u^{[1]}(1) + h u(1) = 0,
        \qquad H,h \in \bR\cup\{\infty\},
\end{equation}
the operator $T=T(\si,H,h)$ becomes~\cite{Sav,SS} selfadjoint, bounded
below, and possesses a simple discrete spectrum accumulating at $+\infty$.
(If $H=\infty$ or $h=\infty$, then the corresponding boundary condition is
understood as a Dirichlet one.)

We observe that among the singular potentials included in this scheme there
are, e.g., the Dirac $\delta$-potentials and Coulomb $1/x$-like potentials
that have been widely used in quantum mechanics to model particle
interactions of various types (see, e.g., \cite{AGHH,AK}). It is also well
known that singular potentials of this kind usually do not produce a single
Sturm-Liouville operator. Still for many reasons $T(\si,H,h)$ can be
regarded as a natural operator associated with differential
expression~\eqref{eq:de} and boundary conditions~\eqref{eq:bc}. For
instance, $l_\si(u) = - u'' + qu$ in the sense of distributions and hence
for regular (i.e., locally integrable) potentials~$q$ the above definition
of $T(\si,H,h)$ coincides with the classical one. Also $T(\si,H,h)$ depends
continuously in the uniform resolvent sense on $\si\in\cH$, which allows us
to regard $T(\si,H,h)$ for singular $q=\si'$ as a limit of regular
Sturm-Liouville operators.

Sturm-Liouville and Schr\"odinger operators with distributional potentials
from the space~$W^{-1}_2(0,1)$ and $W^{-1}_{2,unif}(\bR)$ respectively have
been shown to possess many properties similar to those for operators with
regular potentials, see, e.~g., \cite{HMper,HMgor,Sav,SS}. In particular,
just as in the regular case, the potential $q=\si'$ and boundary
conditions~\eqref{eq:bc} of the Sturm-Liouville operator $T(\si,H,h)$ can
be recovered via the corresponding spectral data, the sequences of
eigenvalues and so-called norming constants~\cite{HMinv}.

The main aim of the present paper is to treat the following inverse
spectral problem. Suppose that $\{\la_n^2\}$ is the spectrum of the
operator~$T(\si,H,h_1)$ and $\{\mu_n^2\}$ is the spectrum of the
operator~$T(\si,H,h_2)$ with $h_2\ne h_1$. Is it possible to recover $\si$,
$H$, $h_1$, and $h_2$ from the two given spectra? More generally, the task
is to find necessary and sufficient conditions for two sequences
$(\la_n^2)$ and $(\mu_n^2)$ so that they could serve as spectra of
Sturm-Liouville operators on $(0,1)$ with some potential $q=\si'\in
W^{-1}_2(0,1)$ and boundary conditions~\eqref{eq:bc} for two different
values of $h$ and, then, to present an algorithm recovering these two
operators (i.e., an algorithm determining~$\si$ and the boundary
conditions).

For the case of a regular (i.e., locally integrable) potential~$q$ the
above problem was treated by \textsc{Levitan and Gasymov}~\cite{Le1,LG}
when $h_1$ and $h_2$ are finite and by
\textsc{Mar\-chen\-ko}~\cite[Ch.~3.4]{Ma} when one of $h_1,h_2$ is
infinite. Earlier, \textsc{Borg}~\cite{Bo} proved that two such spectra
determine the regular potential uniquely.

Observe that $T(\si+h',H-h',h+h')=T(\si,H,h)$ for any $h'\in\bC$, so that
it is impossible to recover $\si,H,h_1$, and $h_2$ uniquely; however, we
shall show that the potential $q=\si'$ and the very operators
$T(\si,H,h_1)$ and $T(\si,H,h_2)$ are determined uniquely by the two
spectra. The reconstruction procedure consists in reduction of the problem
to recovering the potential $q=\si'$ and the boundary conditions based on
the spectrum and the sequence of so-called norming constants. The latter
problem has been completely solved in our paper~\cite{HMinv}, and this
allows us to give a complete description of the set of spectral data and to
develop a reconstruction algorithm for the inverse spectral problem under
consideration.

The organization of the paper is as follows. In the next two sections we
restrict ourselves to the case where $H=\infty$ and $h_1,h_2\in\bR$. In
Section~\ref{sec:asymp} refined eigenvalue asymptotics is found and some
other necessary conditions on spectral data are established. These results
are then used in Section~\ref{sec:reduction} to determine the set of
norming constants and completely solve the inverse spectral problem. The
case where one of $h_1$ and $h_2$ is infinite is treated in
Section~\ref{sec:DN}, and in Section~\ref{sec:N} we comment on the changes
to be made if $H$ is finite. Finally, Appendix~\ref{sec:RB} contains some
facts about Riesz bases of sines and cosines in $L_2(0,1)$ that are
frequently used throughout the paper.


\section{Spectral asymptotics}\label{sec:asymp}

In this section (and until Section~\ref{sec:DN}) we shall consider the case
of the Dirichlet boundary condition at the point $x=0$  and the boundary
conditions of the third type at the point $x=1$ (i.~e., to the case where
$H=\infty$ and $h_1,h_2\in\bR$). We start with the precise definition of
the Sturm-Liouville operators under consideration.

Suppose that $q$ is a real-valued distribution from the class
$W^{-1}_2(0,1)$ and $\si\in\cH$ is any of its real-valued distributional
primitives. For $h\in\bR$ we denote by $T_{\si,h}=T(\si,\infty,h)$ the
operator in $\cH$ given by
\[
    T_{\si,h} u = l_\si(u):= - (u^{[1]})' - \si u'
\]
and the boundary conditions
\[
    u(0) = 0, \qquad u^{[1]}(1) + h u(1) = 0, \qquad h\in\bR
\]
(we recall that $u^{[1]}$ stands for the quasi-derivative~$u'-\si u$
of~$u$). More precisely, the domain of $T_{\si,h}$ equals
\[
    \fD(T_{\si,h}) = \{u\in W^{1}_1[0,1] \mid u^{[1]} \in W^1_1[0,1],\
        l_\si (u) \in \cH,\ u(0)=0,\ u^{[1]}(1) + h u(1) =0\}.
\]
It follows from~\cite{SS} that so defined $T_{\si,h}$ is a selfadjoint
operator with simple discrete spectrum accumulating at $+\infty$. The
results of~\cite{SS1} also imply that the eigenvalues $\la_n^2$ of
$T_{\si,h}$, when ordered increasingly, obey the asymptotic relation $\la_n
= \pi (n-1/2) + \wt\la_n$ with an $\ell_2$-sequence $(\wt\la_n)$. However,
we shall need a more precise form of $\wt\la_n$, and shall derive it next.

To start with, we introduce an operator $T_\si$ defined by $T_\si u =
l_\si(u)$ on the domain
\[
    \fD(T_{\si}) = \{u\in W^{1}_1[0,1] \mid u^{[1]} \in W^1_1[0,1],\
        l_\si (u) \in \cH,\ f(0)=0\}.
\]
In other words, $T_\si$ is a one-dimensional extension of $T_{\si,h}$
discarding the boundary condition at the point $x=1$. We invoke now the
following results of~\cite{HMtr} on the operator~$T_\si$. Just as in the
classical situation with regular potential $q$ the operators~$T_{\pm\si}$
turn out to be similar to the potential-free operator $T_0$; the similarity
is performed by the \emph{transformation operators} $I + K_\si^{\pm}$,
where $K_\si^{\pm}$ are integral Volterra operators of Hilbert-Schmidt
class given by
\[
    K_\si^{\pm} u(x) = \int_0^x k_\si^\pm(x,t) u(t)\, dt.
\]
It follows that any function $u$ from $\fD(T_\si)$ has the form
 $u = (I+K_\si^+)v$ for some function $v \in \fD(T_0)$ (i.e., for some
$v\in W_2^2[0,1]$ satisfying the boundary condition $v(0)=0$) and that
  $T_\si (I+K_\si^+)v = -(I+K_\si^+)v''$.
Moreover, there exists a Hilbert-Schmidt kernel $r_\si(x,t)$ such that
\[
    [(I + K_\si^+)v]^{[1]}(x) = [(I+K_\si^-)v'](x) +
        \int_0^x r_\si(x,t) v(t)\, dt
\]
for every $v \in \fD(T_0)$. Also, for every fixed $x\in[0,1]$ the functions
$k_\si^\pm(x,\cdot)$ and $r_\si(x,\cdot)$ belong to $\cH$.

Due to the above similarity of $T_\si$ and $T_0$, for any nonzero
$\la\in\bC$ the function
\[
    u(x,\la) := \sin \la x + \int_0^x k_\si^+ (x,t) \sin \la t\, dt
\]
is an eigenfunction of the operator $T_\si$ corresponding to the eigenvalue
$\la^2\in\bC$. If in addition $u(x,\la)$ satisfies the boundary condition
$u^{[1]}(1) + h u(1) = 0$, then $u(x,\la)$ is also an eigenfunction of the
operator $T_{\si,h}$ corresponding to the eigenvalue~$\la^2$. Therefore the
nonzero spectrum of~$T_{\si,h}$ consists of the squared nonzero
$\la$-solutions of the equation
\begin{equation}\label{eq:zeros}
 \begin{aligned}
     \cos\la + \int_0^1 k_\si^-(1,t)\cos\la t \,dt
     & + \frac{1}{\la}\int_0^1 r_\si(1,t)\sin\la t \,dt
    \\ & + \frac{h}{\la} \sin\la +
        \frac{h}{\la}\int_0^1 k_\si^+(1,t)\sin\la t \,dt = 0.
 \end{aligned}
\end{equation}
Recall that the operator $T_{\si,h}$ is bounded below and hence after
addition to $q$ a suitable constant $T_{\si,h}$ becomes positive. Therefore
without loss of generality we may assume that all zeros of
equation~\eqref{eq:zeros} are real. Observe also that due to the symmetry
we may consider only the positive zeros. The asymptotics of these zeros is
given in the following theorem.

\begin{theorem}\label{thm:spas}
Suppose that $\si\in\cH$ and $h\in\bR$ are such that the operator
$T_{\si,h}$ is positive and denote by $\la_1^2<\la_2^2<\cdots$ the
eigenvalues of $T_{\si,h}$. Then the numbers $\la_n$ satisfy the relation
\[
    \la_n = \pi (n-1/2) + \frac{h}{\pi n} + \de_n(h) + \ga_n,
\]
in which the sequences $\bigl(\de_n(h)\bigr)_{n=1}^\infty$ and
$(\ga_n)_{n=1}^\infty$ belong to $\ell_1$ and $\ell_2$ respectively with
$\ga_n$ independent of $h$.
\end{theorem}

\begin{proof}
Recall that $\la_n \to \infty$ as $n\to\infty$ and that the functions
$k_\si^\pm(1,t)$ and $r_\si(1,t)$ belong to $\cH$.
Therefore~\eqref{eq:zeros} and the Riemann lemma imply that
\[
    \cos\la_n = -\int_0^1 k_\si^-(1,t)\cos\la_n t \,dt + O(1/\la_n) = o(1)
\]
as $n \to \infty$. By the standard Rouch\'e-type arguments (see details,
e.g., in~\cite[Ch.~1.3]{Ma}) we conclude that $\la_n = \pi (n-1/2) +
\wt\la_n$ with $\wt\la_n \to 0$ as $n\to\infty$. As a result (see
Appendix~\ref{sec:RB}), the system
    $\{\sin\la_n x\}$
turns out to be a Riesz basis of~$\cH$ and the sequences
$(c_n)_{n=1}^\infty$ and $(c'_n)_{n=1}^\infty$ with
\[
    c_n   := \int_0^1 k_\si^+(1,t)\sin\la_nt\,dt, \quad
    c'_n  := \int_0^1 r_\si(1,t)\sin\la_nt\,dt
\]
belong to $\ell_2$. Relation~\eqref{eq:zeros} now implies that
\begin{align*}
    \sin\wt\la_n &= (-1)^{n+1} \int_0^1 k_\si^-(1,t)\cos\la_nt\,dt
        + h \frac{\cos\wt\la_n}{\la_n} + (-1)^{n+1}\frac{hc_n+c'_n}{\la_n} =\\
        &= \ga_n + \frac{h}{\pi n} + \de'_n+\de''_n,
\end{align*}
where
\begin{align}\notag
    \ga_n   &:= (-1)^{n+1} \int_0^1 k_\si^-(1,t)\cos \pi(n-1/2)t\,dt, \\
    \de'_n  &:= (-1)^{n+1} \int_0^1 k_\si^-(1,t)
               [\cos \la_nt -\cos\pi(n-1/2)t]\,dt, \label{eq:rn}\\
    \de''_n &:= \frac{h\cos\wt\la_n}{\la_n} - \frac{h}{\pi n}
                + (-1)^{n+1}\frac{hc_n+c'_n}{\la_n}. \notag
\end{align}

Next we observe that the systems $\{\cos\pi(n-1/2)t\}$ and $\{\cos\la_nt\}$
form Riesz bases of~$\cH$, so that $(\ga_n)$ and $(\de'_n)$ belong to
$\ell_2$. It follows that $(\sin\wt\la_n) \in \ell_2$ and therefore
$(\wt\la_n) \in \ell_2$. Lemma~\ref{lem:fourierdiff} now implies that the
numbers $\de'_n$ form an~$\ell_1$-sequence, and simple arguments justify
the same statement about $\de''_n$.

Thus we have shown that
\[
    \sin\wt\la_n = \ga_n + \frac{h}{\pi n} + \de'_n+ \de''_n
\]
with the sequences $(\ga_n)$ from $\ell_2$ and $(\de'_n)$ and $(\de''_n)$
from $\ell_1$. Applying $\arcsin$ to both parts of this equality, we arrive
at
\[
    \wt\la_n = \ga_n + \frac{h}{\pi n} + \de_n
\]
with some $\ell_1$-sequence $(\de_n)$, and the theorem is proved.
\end{proof}

It was shown in \cite{HMinv} that any sequence $(\la^2_n)$ of pairwise
distinct positive numbers satisfying the relation $\la_n = \pi(n-1/2)+
\wt\la_n$ with an $\ell_2$-sequence $(\wt\la_n)$ is the spectrum of some
Sturm-Liouville operator $T_{\si,h}$ with potential $q=\si'$ from
$W^{-1}_2(0,1)$, the Dirichlet boundary condition at $x=0$, and third type
boundary condition~\eqref{eq:bc} at $x=1$ with some $h\in\bR$. Suppose that
$(\mu_n^2)$ is the spectrum of a Sturm-Liouville operator with the same
potential $q$, the Dirichlet boundary condition at $x=0$, and third type
boundary condition~\eqref{eq:bc} at $x=1$ with a different $h\in\bR$. We
ask whether there exists any relation between the two spectra.

The first restriction, just as in the regular case, is that the spectra
should interlace.

\begin{lemma}\label{lem:interlace}
Suppose that $\si\in\cH$ is real valued and that sequences $\{\la_n^2\}$
and $\{\mu_n^2\}$ are the spectra of Sturm-Liouville operators
$T_{\si,h_1}$ and $T_{\si,h_2}$ with distinct $h_1,h_2\in\ov\bR$. Then
these sequences interlace.
\end{lemma}

\begin{proof}
Denote by $u = u(x,\la)$ a solution of equation $l_\si(u) = \la u$
satisfying the boundary condition $u(0)=0$. We recall that $u$ being a
solution of $l_\si(u) = \la u$ means that
\begin{equation}\label{eq:syst}
    \frac{d}{dx} \bn{u^{[1]}}{u} = \begin{pmatrix}
         -\si & -\si^2 - \la \\ 1 & \si
    \end{pmatrix} \bn{u^{[1]}}{u}
\end{equation}
and hence $u$ enjoys the standard uniqueness properties of solutions to
first order differential systems with regular (i.e. locally integrable)
coefficients. The numbers~$\la_n$ and $\mu_n$ are then solutions of the
equations $u^{[1]}(1,\la) + h_1 u(1,\la) = 0$ and $u^{[1]}(1,\la) + h_2
u(1,\la) = 0$ respectively.

We introduce a continuous function $\th(x,\la)$ through $\cot\th (x,\la)
\equiv {u^{[1]}}/{u}$. After differentiating both sides of this identity in
$x$ and using~\eqref{eq:syst} we get
\[
     -\frac{\th'}{\sin^2{\th}} = \frac{-\si u'u - u^{[1]}u'}{u^2}
                 - \la = -(\cot\th + \si)^2  - \la,
\]
or
\[
    \th' = \la \sin^2\th + (\cos\th + \si\sin\th)^2.
\]
It follows from~\cite[Ch.~8.4]{At} or \cite[proof of Theorem~XI.3.1]{Ha}
that the function $\th(1,\la)$ is strictly increasing in $\la$; hence
solutions of the equations $\cot\th(1,\la)=-h_1$ and $\cot\th(1,\la)=-h_2$
interlace, and the proof is complete.
\end{proof}

The next restriction concerns asymptotics of $\la_n$ and $\mu_n$.

\begin{lemma}\label{lem:diff}
Suppose that $(\la_n^2)$ and $(\mu_n^2)$ are spectra of Sturm-Liouville
operators $T_{\si,h_1}$ and $T_{\si,h_2}$  with real-valued $\si\in\cH$ and
$h_1,h_2\in\bR$. Then there exists an $\ell_2$-sequence $(\nu_n)$ such that
\begin{equation}\label{eq:diff}
    \la_n - \mu_n =  \frac{h_1-h_2}{\pi n} + \frac{\nu_n}{n}.
\end{equation}
\end{lemma}

\begin{proof}
Put $\wt\la_n:= \la_n - \pi(n-1/2)$ and $\wt\mu_n:= \mu_n - \pi(n-1/2)$ so
that $\la_n - \mu_n = \wt\la_n - \wt\mu_n$. Recall that
\[
    \sin\wt\la_n - \sin\wt\mu_n =
      \frac{h_1-h_2}{\pi n} + \de'_n(h_1) - \de'_n(h_2)+
      \de''_n(h_1) - \de''_n(h_2)
\]
with the quantities $\de'_n(h),\de''_n(h)$ of~\eqref{eq:rn}. Relation
\[
   \de'_n(h_1) - \de'_n(h_2) = (-1)^{n+1}
        \int_0^1 k_\si^-(1,t)[\cos \la_nt -\cos\mu_nt]\,dt
\]
and Lemma~\ref{lem:fourierdiff} show that
 \(
    \de'_n(h_1) - \de'_n(h_2) = (\wt\la_n-\wt\mu_n) \nu'_n
 \)
for some $\ell_2$-sequence $(\nu'_n)$. Also there exists an
$\ell_2$-sequence $(\nu''_n)$ such that $\de''_n(h_1)-\de''_n(h_2) =
\nu''_n/n$. Finally, using the relations
\[
    \sin\wt\la_n - \sin\wt\mu_n = (\wt\la_n - \wt\mu_n)\cos\wt\nu_n,
\]
where $\wt\nu_n$ are points between $\wt\la_n$ and $\wt\mu_n$ (so that
$(1-\cos\wt\nu_n)$ is an $\ell_2$-sequence), and combining the above
relations we arrive at
\[
    (\wt\la_n - \wt\mu_n)(\cos\wt\nu_n - \nu'_n) =
       \frac{h_1-h_2}{\pi n} + \frac{\nu''_n}{n},
\]
which implies~\eqref{eq:diff}.
\end{proof}

\begin{corollary}\label{cor:diff}
Under the above assumptions,
 $\lim_{n\to\infty}(\la_n^2 - \mu_n^2) = 2(h_1-h_2)$.
\end{corollary}


\section{Reduction to the inverse spectral problem by one spectrum
and norming constants}\label{sec:reduction}

Suppose that $\si\in\cH$ is real-valued and that $\{\la^2_n\}$ and
$\{\mu^2_n\}$ are eigenvalues of the operators $T_{\si,h_1}$ and
$T_{\si,h_2}$ respectively introduced in the previous section. We shall
show how the problem of recovering $\si$, $h_1$, and $h_2$ based on these
spectra can be reduced to the problem of recovering $\si$ and $h_1$ based
on the spectrum $\{\la^2_n\}$ of $T_{\si,h_1}$ and the so-called norming
constants~$\{\al_n\}$ defined below. This latter problem for the class of
Sturm-Liouville operators with singular potentials from $W^{-1}_2(0,1)$ is
completely solved in~\cite{HMinv}; see also~\cite{Le,Ma} for the regular
case of integrable potentials.

Denote by $u_1(x,\la)$ and $u_2(x,\la)$ solutions of the equation
 $l_\si(u) = \la^2 u$
satisfying the initial conditions $u_1(1,\la)=u_2(1,\la)=1$ and
$u_j^{[1]}(1,\la) = - h_j$, $j=1,2$. Then $u_1(x,\la_n)$ are the
eigenfunctions of the operator $T_{\si,h_1}$ corresponding to the
eigenvalues $\la_n^2$, and we put
\[
    \al_n := 2 \int_0^1 |u_1(x,\la_n)|^2\,dx.
\]
Our next aim is to show that $\al_n$ can be expressed in terms of the
spectra $\{\la^2_k,\mu^2_k\}$ only.

Set $\Phi_j(\la):=u_j(0,\la)$, $j=1,2$; then zeros of $\Phi_1$ and $\Phi_2$
are precisely the numbers $\pm\la_n$ and $\pm\mu_n$ respectively. We show
that $\Phi_1$ and $\Phi_2$ uniquely determine $\al_n$ and are uniquely
determined by their zeros.

\begin{lemma}\label{lem:al}
The norming constants $\al_n$ satisfy the following equality:
\begin{equation}\label{eq:al}
    \al_n = \frac{h_1-h_2}{\la_n}\frac{\Phi'_1(\la_n)}{\Phi_2(\la_n)}.
\end{equation}
\end{lemma}

\begin{proof}
The proof is rather standard and the details can be found, e.g.,
in~\cite{Le}, so we shall only sketch its main points here.

Put $m(\la):= -\Phi_2(\la)/\Phi_1(\la)$; then the function
$f(x,\la):=u_2(x,\la)+m(\la)u_1(x,\la)$ for all $\la\in\bC$ not in the
spectrum of~$T_{\si,h_1}$ satisfies the equation $l_\si(f) = \la^2 f$ and
the boundary condition $f(0,\la)=0$. In other words, $f(x,\la)$ is an
eigenfunction of the operator $T_\si$ corresponding to the eigenvalue
$\la^2$. Since $u_1(x,\la_n)$ is an eigenfunction of the operator~$T_\si$
corresponding to the eigenvalue $\la_n^2$, we have
\begin{align*}
    (\la^2-\la^2_n)\int_0^1 f(x,\la) u_1(x,\la_n)\,dx &=
        (T_\si f,u_1) - (f,T_\si u_1)\\
        &= - f^{[1]}(1,\la)\ov{u_1(1,\la)} + f(1,\la)\ov{u^{[1]}(1,\la_n)}
            = h_2 - h_1.
\end{align*}
On the other hand,
\begin{align*}
    (\la^2-\la^2_n)\int_0^1 f(x,\la) u_1(x,\la_n)\,dx &=
        (\la^2-\la_n^2)\int_0^1 u_2(x,\la)u_1(x,\la_n)\,dx \\
        &- (\la^2-\la_n^2)\frac{\Phi_2(\la)}{\Phi_1(\la)}
                \int_0^1 u_1(x,\la)u_1(x,\la_n)\,dx,
\end{align*}
and after combining the two relations and letting $\la\to\la_n$ we arrive
at
\[
    \al_n = \frac{h_1-h_2}{\la_n}\frac{\Phi'_1(\la_n)}{\Phi_2(\la_n)}.
\]
The proof is complete.
\end{proof}

The following statement has appeared in many variants in numerous sources,
but we include its proof here for the sake of completeness.

\begin{lemma}\label{lem:repr}
In order that a function $f(\la)$ admit the representation
\begin{equation}\label{eq:int}
    f(\la) = \cos\la + \int_0^1 g(t)\cos\la t\,dt
\end{equation}
with an $\cH$-function $g$, it is necessary and sufficient that
\begin{equation}\label{eq:prod}
    f(\la) = \prod_{k=1}^\infty \frac{(f_k^2 - \la^2)}{\pi^2(k-1/2)^2},
\end{equation}
where $f_k := \pi(k-1/2) + \wt f_k$ and $(\wt f_k)$ is an
$\ell_2$-sequence.
\end{lemma}

\begin{proof}
\emph{Necessity.} The function $f$ of~\eqref{eq:int} is an even entire
function of order~$1$, and the standard Rouch\'e-type arguments (see the
analysis of Section~\ref{sec:asymp}) show that zeros~$\pm f_k$ of~$f$ have
the asymptotics $f_k = \pi (k-1/2) + \wt f_k$ with $(\wt f_k)\in\ell_2$. On
the other hand, up to a scalar factor $C$, the function~$f$ is recovered
from its zeros as
\begin{equation}\label{eq:c1}
    f(\la) = C\prod_{k=1}^\infty \Bigl(1 - \frac{\la^2}{f_k^2} \Bigr)
        = C_1 \prod_{k=1}^\infty \frac{(f_k^2 - \la^2)}{\pi^2(k-1/2)^2},
\end{equation}
where
\[
    C_1 = C\prod_{k=1}^\infty \frac{\pi^2(k-1/2)^2}{f_k^2}.
\]
To determine $C_1$, we observe that, in view of~\eqref{eq:int},
 $\lim_{\nu\to\infty}f(i\nu)/\cos i\nu=1$
and that
\begin{equation}\label{eq:cos}
    \cos\la = \prod_{k=1}^{\infty}\frac{(\pi^2(k-1/2)^2 - \la^2)}
            {\pi^2(k-1/2)^2};
\end{equation}
comparing now~\eqref{eq:c1} and \eqref{eq:cos}, we conclude that $C_1=1$
and necessity is justified.

\emph{Sufficiency.} Suppose now that $f$ has a representation of the
form~\eqref{eq:prod}, in which $f_k = \pi(k-1/2)+ \wt f_k$ with
 $(\wt f_k)\in\ell_2$. We assume first that the zeros~$f_k$ are pairwise
distinct. Then due to the asymptotics of $f_k$ the system
 $(\cos f_kx)_{k=1}^\infty$ is a Riesz basis of~$\cH$
 (see Lemma~\ref{lem:RB}) and the numbers
 $c_k:=-\cos f_k=(-1)^{k}\sin \wt f_k$ define an $\ell_2$-sequence.
Therefore there exists a unique function $g\in\cH$ with Fourier
coefficients $c_k$ in the Riesz basis $(\cos f_kx)_{k=1}^\infty$. Now the
function
\begin{equation}\label{eq:wtf}
    \wt f(\la):= \cos \la + \int_0^1 g(x)\cos\la x\,dx
\end{equation}
is an even entire function of order~$1$ and has zeros at the points
 $\pm f_k$. It follows that $\wt f$ and $f$ differ by a scalar factor,
which as before is shown to equal~$1$.

Due to the asymptotics of the zeros~$f_k$ only finitely many of them can
repeat. If, e.g., $f_n = f_{n+1} = \cdots = f_{n+p-1}$ is a zero of $f$ of
multiplicity~$p=p(n)$, then we include to the above system the functions
$\cos f_nx, x\sin f_nx, \cdots, x^{p-1}\sin(\pi p/2 - f_nx)$. After this
modification has been done for every multiple root, we again get a Riesz
basis of $\cH$ and can find a function $g\in\cH$ with Fourier coefficients
$c_n:=-\cos f_n, c_{n+1} = -\sin f_n, \cdots,
    c_{n+p(n)}:= -\sin(\pi p(n)/2- f_nx)$.
Then $\pm f_n$ becomes a zero of the function $\wt f$~in~\eqref{eq:wtf} of
multiplicity $p(n)$, and the proof is completed as above.
\end{proof}

\begin{lemma}\label{lem:phi}
The following equalities hold:
\begin{equation}\label{eq:phi}
 \begin{gathered}
    \Phi_1(\la) = \prod\pi^{-2}(k-1/2)^{-2}(\la^2_k-\la^2), \\
    \Phi_2(\la) = \prod\pi^{-2}(k-1/2)^{-2}(\mu^2_k-\la^2).
 \end{gathered}
\end{equation}
\end{lemma}

\begin{proof}
We recall~\cite{HMtr} that the function $u_1$ has a representation via a
transformation operator connected with the point $x=1$ of the form
\[
    u_1(x,\la) = \cos\la(1-x) + \int_x^1k(x,t)\cos\la (1-t)\,dt,
\]
which implies that
\[
    \Phi_1(\la) = \cos\la + \int_0^1 k_1(t)\cos\la t\,dt,
\]
where $k_1(t):= k(0,1-t)$ is a function from $\cH$. The claim now follows
from Lemma~\ref{lem:repr}. The representation for $\Phi_2$ is derived
analogously, and the lemma is proved.
\end{proof}

Now, given two spectra $(\la^2_n)$ and $(\mu^2_n)$ of two Sturm-Liouville
operators $T_{\si,h_1}$ and $T_{\si,h_2}$ respectively, with unknown
$\si\in\cH$ and $h_1,h_2\in\bR$, we proceed as follows. First, we identify
$h_1-h_2$ as $ \lim_{n\to\infty} (\la^2_n-\mu^2_n)/2$ (see
Corollary~\ref{cor:diff}), then construct the functions $\Phi_1$ and
$\Phi_2$ of~\eqref{eq:phi} and determine the norming coefficients $(\al_n)$
via~\eqref{eq:al}. The spectral data $\{(\la^2_n), (\al_n)\}$ determine
$\si$ and $h_1$ up to an additive constant by means of the algorithm
of~\cite{HMinv}. This gives the required function $\si$ and two numbers
$h_1$ and $h_2$ up to an additive constant, and the reconstruction is
complete.

The second part of the inverse spectral problem is to identify those pairs
of sequences $(\la_n^2)$ and $(\mu^2_n)$ that are spectra of
Sturm-Liouville operators $T_{\si,h_1}$ and $T_{\si,h_2}$ for some
real-valued $\si\in\cH$ and some real $h_1,h_2$, $h_1\ne h_2$. We
established in Section~\ref{sec:asymp} the necessary conditions on the two
spectra given by Theorem~\ref{thm:spas} and Lemmata~\ref{lem:interlace} and
\ref{lem:diff}; it turns out that these conditions are sufficient as well.
Namely, the following statement holds true.

\begin{theorem}\label{thm:suff}
Suppose that sequences $(\la^2_n)$ and $(\mu^2_n)$ of positive pairwise
distinct numbers satisfy the following assumptions:
\begin{itemize}
\item [(1)] the sequences $(\la^2_n)$ and $(\mu^2_n)$ interlace;
\item [(2)] $\la_n = \pi(n-1/2) + \wt\la_n$ and
            $\mu_n = \pi(n-1/2) + \wt\mu_n$ with some
            $\ell_2$-sequences~$(\wt\la_n)$ and $(\wt\mu_n)$;
\item [(3)] there exist a real number $h$ and an $\ell_2$-sequence $(\nu_n)$
            such that $\la_n-\mu_n=\frac{h}{\pi n}+\frac{\nu_n}{n}$.
\end{itemize}
Then there exist a function $\si\in\cH$ and two real numbers $h_1$ and
$h_2$ such that $(\la^2_n)$ and $(\mu^2_n)$ are the spectra of the
Sturm-Liouville operators $T_{\si,h_1}$ and $T_{\si,h_2}$ respectively.
Moreover, the operators $T_{\si,h_1}$ and $T_{\si,h_2}$ are recovered
uniquely.
\end{theorem}

In order to prove the theorem we have to show first that the
numbers~$\al_n$ constructed for the two sequences through
formula~\eqref{eq:al}, with the functions~$\Phi_1,\Phi_2$ of~\eqref{eq:phi}
and $h_1-h_2:=h$ with $h$ of item~(3), are positive and obey the
asymptotics $\al_n = 1 + \wt\al_n$ for some $\ell_2$-sequence $(\wt\al_n)$.
Then the algorithm of~\cite{HMinv} uses $(\la^2_n)$ and $(\al_n)$ to
determine (unique up to an additive constant) function $\si\in\cH$ and
$h_1\in\bR$ such that $\{\la_n^2\}$ is the spectrum of the Sturm-Liouville
operator~$T_{\si,h_1}$ and $\al_n$ are the corresponding norming constants.
The second and final step is to verify that $\{\mu^2_n\}$ is the spectrum
of the Sturm-Liouville operator $T_{\si,h_2}$ with $h_2:= h_1 - h$, where
$h$ is the number of item~(3). These two steps are performed in the
following two lemmata.

\begin{lemma}\label{lem:alas}
Assume (1)--(3) of Theorem~\ref{thm:suff}. Then the numbers $\al_n$
constructed through relation~\eqref{eq:al} with $\Phi_1,\Phi_2$
of~\eqref{eq:phi} and $h_1-h_2:=h$ with $h$ of item~(3) are all positive
and obey the asymptotics
 \(
    \al_n = 1 + \wt\al_n,
 \)
where $(\wt\al_n)$ is an $\ell_2$-sequence.
\end{lemma}

\begin{proof}
By Lemma~\ref{lem:repr} there exist $\cH$-functions $f_1$, $f_2$ such that
the functions $\Phi_j$, $j=1,2$, of~\eqref{eq:phi} admit the representation
\[
    \Phi_j(\la) = \cos\la + \int_0^1 f_j(t) \cos\la t\, dt.
\]
It follows that
\begin{align*}
    \Phi_1'(\la_n) &= -\sin\la_n - \int_0^1 tf_1(t) \sin\la_nt\,dt\\
        &= (-1)^n \cos\wt\la_n - \int_0^1 tf_1(t) \sin\la_nt\,dt
        = (-1)^n + \hat\la_n,
\end{align*}
where $(\hat\la_n)_{n=1}^\infty \in \ell_2$, and by similar arguments
\begin{align*}
    \Phi_2(\la_n) &= (\la_n-\mu_n)\Phi_2'(\mu_n) + O(|\la_n-\mu_n|^2)\\
        & = (\la_n-\mu_n)[(-1)^n + \hat\mu_n]
\end{align*}
for some $\ell_2$-sequence $(\hat\mu_n)$. Next, assumptions (2) and (3)
easily imply the relation
\[
    \frac{h}{\la_n(\la_n-\mu_n)} = 1 + \hat\nu_n,
\]
for some $\ell_2$-sequence~$(\hat\nu_n)$. Therefore the numbers
\[
    \al_n:= \frac{h}{\la_n}\frac{\Phi'_1(\la_n)}{\Phi_2(\la_n)}
\]
obey the required asymptotics.

Finally, the interlacing property of the two sequences implies that all
$\al_n$ are of the same sign, and thus are all positive in view of the
asymptotics established. The proof is complete.
\end{proof}

With the above-stated properties of the sequences $(\la_n^2)$ and $(\al_n)$
we can employ the result of~\cite{HMinv} that guarantees existence of a
unique Sturm-Liouville operator $T_{\si,h_1}$ with a real-valued
function~$\si\in\cH$ and a real number $h_1$ such that $\{\la^2_n\}$
coincides with the spectrum of $T_{\si,h_1}$ and $\al_n$ are the
corresponding norming constants. Now we put $h_2:=h_1-h$ with $h$ of
assumption~(3) and expect $\{\mu^2_n\}$ to be the spectrum of
$T_{\si,h_2}$.

\begin{lemma}\label{lem:mu}
The spectrum of the Sturm-Liouville operator $T_{\si,h_2}$ with the above
$\si\in\cH$ and $h_2\in\bR$ coincides with $(\mu^2_n)$.
\end{lemma}

\begin{proof}
Denote by $v_1(x,\la)$ a solution of the equation $l_\si(v) = \la^2v$
satisfying the boundary conditions $v(1)=1, v^{[1]}(1) = -h_1$. Then
$v_1(x,\la_n)$ is an eigenfunction of the operator $T_{\si,h_1}$
corresponding to the eigenvalue $\la^2_n$, and by construction
\[
    2\int_0^1 |v_1(x,\la_n)|^2\, dx = \al_n :=
        \frac{h}{\la_n}\frac{\Phi'_1(\la_n)}{\Phi_2(\la_n)}
\]
with the functions $\Phi_1,\Phi_2$ of~\eqref{eq:phi}.

Denote by $(\nu^2_n)$ the spectrum of $T_{\si,h_2}$ and put
\[
    \wt\Phi_2(\la):=\prod_{k=1}^\infty \pi^{-2}(k-1/2)^{-2}(\nu^2_n-\la^2).
\]
By Lemma~\ref{lem:al} we have
\[
    \al_n =\frac{h}{\la_n}\frac{\Phi'_1(\la_n)}{\wt\Phi_2(\la_n)},
\]
whence $\Phi_2(\la_n) = \wt\Phi_2(\la_n)$ for all $n\in\bN$.
Lemma~\ref{lem:repr} implies that
\[
    \Phi_2(\la) - \wt\Phi_2(\la) = \int_0^1 f(t) \cos\la t\,dt
\]
for some $f\in\cH$. Equality $\Phi_2(\la_n) = \wt\Phi_2(\la_n)$ means that
the function $f$ is orthogonal to $\cos\la_nt$, $n\in\bN$. Since the system
$\{\cos\la_nt\}$ forms a Riesz basis of $\cH$, we get $f\equiv0$ and
$\Psi_2 \equiv \wt\Psi_2$. Thus $\nu_n = \mu_n$, and the lemma is proved.
\end{proof}

\section{Reconstruction by Dirichlet and Dirichlet-Neumann
    spectra}\label{sec:DN}

The analysis of the previous two sections does not cover the case where one
of $h_1,h_2$ is infinite. In this case the other number may be taken~$0$
without loss of generality (recall that $T_{\si,h}=T_{\si+h,0}$), i.e., the
boundary conditions under considerations become Dirichlet and
Dirichlet-Neumann ones.

Suppose therefore that $\si\in\cH$ is real valued and that $(\la_n^2)$ and
$(\mu_n^2)$ are spectra of the operators $T_{\si,0}$ and $T_{\si,\infty}$
respectively; without loss of generality we assume that $\la_n$ and $\mu_n$
are positive and strictly increase with $n$. The reconstruction procedure
remains the same as before; namely, we use the two spectra to determine a
sequence of norming constants $(\al_n)$ and then recover $\si$ by the
spectral data $\{(\la_n^2),(\al_n)\}$.

Denote by $u_-(\cdot,\la)$ and $u_+(\cdot,\la)$ solutions of the equation
$l_\si u = \la u$ satisfying the initial conditions $u_-(0,\la)=0,
u^{[1]}_-(0,\la)=\la$ and $u_+(1,\la)=1,u^{[1]}_+(1,\la)=0$ respectively.
Then $u_+(\cdot,\la_n)$ is an eigenfunction of the operator $T_{\si,0}$
corresponding to the eigenvalue~$\la_n^2$ and
\[
    \al_n:= 2 \int_0^1 |u_+(x,\la_n)|^2\,dx
\]
is the corresponding norming constant. We also put
$\Psi_1(\la):=u_+(0,\la)$ and $\Psi_2(\la):=u_-(1,\la)$; then zeros of the
functions $\Psi_1$ and $\Psi_2$ are numbers $\pm\la_n$ and $\pm\mu_n$
respectively. As earlier, the function $\Psi_1$ is uniquely determined by
its zeros through formula~\eqref{eq:phi}. Observe that the Dirichlet
eigenvalues $\mu_n^2$ have asymptotics different from that of $\la_n^2$, so
that~$\Psi_2$ requires slight modification of formula~\eqref{eq:phi}. We
first investigate the asymptotics of~$\mu_n$.

\begin{theorem}[{\cite{HMtr,Sav}}]\label{thm:spas1}
Suppose that $\si\in\cH$ and $\{\mu_n^2\}$, $n\in\bN$, is the spectrum of
the operator $T_{\si,\infty}$. Then the numbers $\mu_n$ satisfy the
relation
\[
    \mu_n = \pi n + \wt\mu_n,
\]
in which the sequence $(\wt\mu_n)$ belongs to $\ell_2$.
\end{theorem}

\begin{proof}
The solution $u_-$ can be represented by means of the transformation
operator~\cite{HMtr} as $u_-(x,\la)=\sin\la x + \int_0^x k_\si^+ (x,t) \sin
\la t\, dt$, where $k_\si^+$ is the kernel of the transformation operator.
Therefore the numbers $\pm\mu_n$ are zeros of the function
\begin{equation}\label{eq:zerosinfty}
    \Psi_2(\la):= \sin\la + \int_0^1 k_\si^+ (1,t) \sin \la t\, dt,
\end{equation}
which is entire of order~$1$. Since the function $k_\si^+(1,\cdot)$ belongs
to~$\cH$~\cite{HMtr}, the required asymptotics of $\mu_n$ is derived in a
standard way (cf.~\cite{Ma} and Section~\ref{sec:asymp}).
\end{proof}

The function $\Psi_2$ is determined by its its zeros in the following way.

\begin{lemma}\label{lem:repr'}
The following equality holds:
\begin{equation}\label{eq:phi'}
      \Psi_2(\la) = \la\prod\pi^{-2}k^{-2}(\mu^2_k-\la^2).
\end{equation}
\end{lemma}

The proof is completely analogous to that of Lemma~\ref{lem:repr} and is
left to the reader.

Now we show how the norming constants $\al_n$ are expressed via $\Psi_1$
and~$\Psi_2$.

\begin{lemma}\label{lem:al'}
The norming constants satisfy the following equality:
\begin{equation}\label{eq:al'}
    \al_n = - \frac{\Psi'_1(\la_n)}{\Psi_2(\la_n)}.
\end{equation}
\end{lemma}

\begin{proof}[Proof (cf.~\cite{GS})]
The Green function $G(x,y,\la^2)$ of the operator~$T_{\si,0}$ (i.~e., the
kernel of the resovlent $(T_{\si,0}-\la^2)^{-1}$) equals
\[
    G(x,y,\la^2) = \sum_{n=1}^{\infty}
        \frac2{\la^2_n - \la^2}\frac{u_+(x,\la_n)u_+(y,\la_n)}{\al_n}
\]
On the other hand, we have
\[
    G(x,y,\la^2) = \frac1{W(\la)} \cdot \left\{
        \begin{aligned} u_-(x,\la) u_+(y,\la),& \qquad
                    0\le x\le y\le1,\\
                           u_+(x,\la)u_-(y,\la),& \qquad
                    0\le y<x\le1,
        \end{aligned} \right.
\]
where $W(\la):=u_+(x,\la)u^{[1]}_-(x,\la)-u_-(x,\la)u^{[1]}_+(x,\la)$ is
the Wronskian of $u_+$ and $u_-$. The value of $W(\la)$ is independent of
$x\in[0,1]$; in particular, taking $x=0$ we find that $W(\la)\equiv \la
u_+(0,\la)=\la \Psi_1(\la)$.

Now we take $x=y=1$ in the above expressions and find that
\[
    \sum_{n=1}^\infty \frac2{\la_n^2-\la^2}\frac{1}{\al_n}
            \equiv \frac{\Psi_2(\la)}{\la\Psi_1(\la)}.
\]
Comparing the residues at the poles $\la=\la_n$, we derive
formula~\eqref{eq:al'}, and the lemma is proved.
\end{proof}

Now the reconstruction procedure is completed as follows: we determine the
sequence of norming constants $(\al_n)$ as explained above and find a
unique operator $T_{\si,0}$ with spectrum $(\la_n^2)$ and norming
constants~$(\al_n)$, see~\cite{HMinv}. This gives the function $\si\in\cH$
and thus the operator $T_{\si,\infty}$.

Now we would like to give an explicit description of the set of all
possible spectra of the operators $T_{\si,0}$ and $T_{\si,\infty}$ when a
real-valued function $\si$ runs through $\cH$, i.~e., to give the necessary
and sufficient conditions on two sequences $(\la_n^2)$ and $(\mu_n^2)$ to
be the spectra of~$T_{\si,0}$ and $T_{\si,\infty}$ with a real-valued
$\si\in\cH$. Necessary conditions are that the eigenvalues $\la_n^2$ and
$\mu_n^2$ should interlace and obey the asymptotics described in
Theorems~\ref{thm:spas} and \ref{thm:spas1}; we shall show that these
conditions are in fact sufficient as well.

\begin{theorem}\label{thm:suffDD}
Suppose that sequences $(\la^2_n)$ and $(\mu^2_n)$ of positive pairwise
distinct numbers satisfy the following assumptions:
\begin{itemize}
\item [(1)] the sequences $(\la^2_n)$ and $(\mu^2_n)$ interlace:
    $0<\la_1^2<\mu_1^2<\la_2^2<\dots$;
\item [(2)] $\la_n = \pi(n-1/2) + \wt\la_n$ and
            $\mu_n = \pi n + \wt\mu_n$ with some
            $\ell_2$-sequences~$(\wt\la_n)$ and $(\wt\mu_n)$.
\end{itemize}
Then there exist a unique function $\si\in\cH$ such that $(\la^2_n)$ and
$(\mu^2_n)$ are the spectra of the Sturm-Liouville operators $T_{\si,0}$
and $T_{\si,\infty}$ respectively.
\end{theorem}

The main ingredients of the proof of this theorem are the same as for
Theorem~\ref{thm:suff}: first, we construct functions $\Psi_1$ and $\Psi_2$
by their zeros and define the sequence $(\al_n)$ as explained in
Lemma~\ref{lem:al'}; then we prove that $\al_n$ are all positive and have
the required asymptotics $\al_n = 1+ \wt\al_n$ for some $\ell_2$-sequence
$(\wt\al_n)$. Using now the reconstruction algorithm of~\cite{HMinv}, we
find a unique function $\si\in\cH$ such that the corresponding
Sturm-Liouville operator $T_{\si,0}$ with potential $q=\si'\in\cH$
possesses the spectral data~$\{(\la_n^2),(\al_n)\}$. Finally, we prove that
$(\mu_n^2)$ is the spectrum of the operator~$T_{\si,\infty}$ with $\si$
just found, and  the reconstruction procedure is complete.

\section{The case of the Neumann boundary condition at $x=0$}\label{sec:N}

The analysis of the previous sections can easily be modified to cover the
case $H=0$, i.e., the Neumann boundary condition $u^{[1]}(0)=0$. As before,
we use the two spectra to determine the sequence of norming constants and
then apply the reconstruction procedure of~\cite{HMinv} to find the
corresponding Sturm-Liouville operators. Also the necessary and sufficient
conditions on the two spectra can be established. We formulate the
corresponding results in the following two theorems.

\begin{theorem}\label{thm:suffNN}
Suppose that  sequences $(\la^2_n)$ and $(\mu^2_n)$ of positive pairwise
distinct numbers satisfy the following assumptions:
\begin{itemize}
\item [(1)] the sequences $(\la^2_n)$ and $(\mu^2_n)$ interlace;
\item [(2)] $\la_n = \pi(n-1) + \wt\la_n$ and
            $\mu_n = \pi(n-1) + \wt\mu_n$ with some
            $\ell_2$-sequences~$(\wt\la_n)$ and $(\wt\mu_n)$;
\item [(3)] there exist a real number $h$ and an $\ell_2$-sequence $(\nu_n)$
            such that $\la_n-\mu_n=\frac{h}{\pi n}+\frac{\nu_n}{n}$.
\end{itemize}
Then there exist a unique function $\si\in\cH$ and real constants $h_1,h_2$
such that $(\la^2_n)$ and $(\mu^2_n)$ are the spectra of the
Sturm-Liouville operators $T_{\si,0,h_1}$ and $T_{\si,0,h_2}$ respectively.

Conversely, the spectra $(\la_n^2)$ and $(\mu_n^2)$ of Sturm-Liouville
operators $T_{\si,0,h_1}$ and $T_{\si,0,h_2}$ with $\si\in\cH$ and
$h_1,h_2\in\bR$ satisfy assumptions~(1)--(3).
\end{theorem}

For the case where one of $h_1,h_2$ is infinite (say, $h_2=\infty$) the
asymptotics of the corresponding spectrum is different; also assumption~(3)
becomes meaningless and should be omitted.

\begin{theorem}\label{thm:suffND}
Suppose that  sequences $(\la^2_n)$ and $(\mu^2_n)$ of positive pairwise
distinct numbers satisfy the following assumptions:
\begin{itemize}
\item [(1)] the sequences $(\la^2_n)$ and $(\mu^2_n)$ interlace;
\item [(2)] $\la_n = \pi(n-1) + \wt\la_n$ and
            $\mu_n = \pi(n-1/2) + \wt\mu_n$ with some
            $\ell_2$-sequences~$(\wt\la_n)$ and $(\wt\mu_n)$.
\end{itemize}
Then there exist a unique function $\si\in\cH$ and a real constant $h_1$
such that $(\la^2_n)$ and $(\mu^2_n)$ are the spectra of the
Sturm-Liouville operators $T_{\si,0,h_1}$ and $T_{\si,0,\infty}$
respectively.

Conversely, the spectra $(\la_n^2)$ and $(\mu_n^2)$ of Sturm-Liouville
operators $T_{\si,0,h_1}$ and $T_{\si,0,\infty}$ with $\si\in\cH$ and
$h_1\in\bR$ satisfy assumptions~(1) and (2).
\end{theorem}

\appendix
\section{Riesz bases}\label{sec:RB}

In this appendix we gather some well known facts about Riesz bases of sines
and cosines (see, e.g., \cite{GK1,HV} and references therein for a detailed
exposition of this topic).

Recall that a sequence $(e_n)_{1}^\infty$ in a Hilbert space $\cH$ is a
Riesz basis if and only if any element $e\in\cH$ has a unique expansion
        \( e = \sum_{n=1}^\infty c_n e_n\)
with $(c_n)\in \ell_2$. If $(e_n)$ is a Riesz basis, then in the above
expansion the \emph{Fourier coefficients} $c_n$ are given by $c_n =
(e,e'_n)$, where $(e'_n)_{1}^\infty$ is a system biorthogonal to $(e_n)$,
i.e., a system which satisfies the equalities $(e_k,e'_n) = \de_{kn}$ for
all $k,n\in\bN$. Moreover, the biorthogonal system $(e'_n)$ is a Riesz
basis of $\cH$ as long as $(e_n)$ is, in which case for any $e\in\cH$ also
the expansion $e=\sum(e,e_n)e'_n$ takes place. In particular, if $(e_n)$ is
a Riesz basis, then for any $e\in\cH$ the sequence $(c'_n)$ with
$c'_n:=(e,e_n)$ belongs to $\ell_2$.

\begin{proposition}[{\cite{HV}}]\label{lem:RB}
Suppose that $\mu_k\to0$ as $k\to\infty$ and that the sequence $\pi
k+\mu_k$ is strictly increasing. Then each of the following systems forms a
Riesz basis of $L_2(0,1)$:
\begin{itemize}
\item [(a)] $\{\sin(\pi kx + \mu_kx)\}_{k=1}^\infty$;
\item [(b)] $\{\sin(\pi [k-1/2]x + \mu_kx)\}_{k=1}^\infty$;
\item [(c)] $\{\cos(\pi kx + \mu_kx)\}_{k=0}^\infty$;
\item [(d)] $\{\cos(\pi [k+1/2]x + \mu_kx)\}_{k=0}^\infty$.
\end{itemize}
\end{proposition}

\begin{lemma}\label{lem:fourierdiff}
Suppose that $(\la_n)$ and $(\mu_n)$ are two sequences of real numbers such
that
\[
    \lim_{n\to\infty} \bigl[\la_n-\pi(n-1/2)\bigr] =
    \lim_{n\to\infty} \bigl[\mu_n-\pi(n-1/2)\bigr] = 0.
\]
and assume that $f\in\cH$. Then there exists an $\ell_2$-sequence $(\nu_n)$
such that
\[
    \int_0^1 f(t) [\cos\la_nt - \cos \mu_nt]\,dt =
            (\la_n-\mu_n)\nu_n + O(|\la_n-\mu_n|^3)
\]
as $n\to\infty$.
\end{lemma}

\begin{proof}
Using the relation
\[
    \cos\la_nt-\cos\mu_nt = -
    2\sin[(\la_n-\mu_n)t/2]\sin[(\la_n+\mu_n)t/2],
\]
we find that
\begin{multline*}
    \int_0^1 f(t) [\cos\la_nt - \cos \mu_nt]\,dt
        = -2\int_0^1f(t)\sin[(\la_n-\mu_n)t/2]\sin[(\la_n+\mu_n)t/2]\,dt\\
        = -(\la_n -\mu_n) \int_0^1 tf(t) \sin[(\la_n+\mu_n)t/2]\,dt +
            O(|\la_n-\mu_n|^3).
\end{multline*}
We put $\nu_n:=-\int_0^1 tf(t) \sin[(\la_n+\mu_n)t/2]\,dt$ and observe that
the sequence $(\nu_n)$ belongs to $\ell_2$ since the system
$\{\sin[(\la_n+\mu_n)t/2]\}$ is Riesz basic in~$\cH$ (at least for all $n$
large enough) and the function $tf(t)$ belongs to $\cH$.
\end{proof}


\end{document}